\documentclass[12pt]{amsart}
\usepackage{amssymb, amsmath, amscd}

\input cyracc.def

\newfont{\rmm}{cmr10 scaled 1000}
\newfont{\bff}{cmbx10  scaled 1000}
\newfont{\itt}{cmsl10 scaled 1000}

\newcounter{fig}

\begin{document}

\renewcommand{\a}{\alpha}
\renewcommand{\b}{\beta}
\newcommand{\g}{\gamma}
\renewcommand{\d}{\delta}
\newcommand{\e}{\epsilon}
\newcommand{\ve}{\varepsilon}
\newcommand{\z}{\zeta}
\renewcommand{\t}{\theta}
\renewcommand{\l}{\lambda}
\renewcommand{\k}{\varkappa}
\newcommand{\m}{\mu}
\newcommand{\n}{\nu}
\renewcommand{\r}{\rho}
\newcommand{\vr}{\varrho}
\newcommand{\s}{\sigma}
\newcommand{\vp}{\varphi}
\renewcommand{\o}{\omega}

\newcommand{\G}{\Gamma}
\newcommand{\D}{\Delta}
\newcommand{\T}{\Theta}
\renewcommand{\L}{\Lambda}
\renewcommand{\P}{\Pi}
\newcommand{\Si}{\Sigma}
\renewcommand{\O}{\Omega}

\renewcommand{\AA}{{\mathcal A}}
\newcommand{\BB}{{\mathcal B}}
\newcommand{\CC}{{\mathcal C}}
\newcommand{\DD}{{\mathcal D}}
\newcommand{\EE}{{\mathcal E}}
\newcommand{\FF}{{\mathcal F}}
\newcommand{\GG}{{\mathcal G}}
\newcommand{\HH}{{\mathcal H}}
\newcommand{\II}{{\mathcal I}}
\newcommand{\JJ}{{\mathcal J}}
\newcommand{\KK}{{\mathcal K}}
\newcommand{\LL}{{\mathcal L}}
\newcommand{\MM}{{\mathcal M}}
\newcommand{\NN}{{\mathcal N}}
\newcommand{\OO}{{\mathcal O}}
\newcommand{\PP}{{\mathcal P}}
\newcommand{\QQ}{{\mathcal Q}}
\newcommand{\RR}{{\mathcal R}}
\renewcommand{\SS}{{\mathcal S}}
\newcommand{\TT}{{\mathcal T}}
\newcommand{\UU}{{\mathcal U}}
\newcommand{\VV}{{\mathcal V}}
\newcommand{\WW}{{\mathcal W}}
\newcommand{\XX}{{\mathcal X}}
\newcommand{\YY}{{\mathcal Y}}
\newcommand{\ZZ}{{\mathcal Z}}

\renewcommand{\aa}{{\mathbb{A}}}
\newcommand{\bb}{{\mathbb{B}}}
\newcommand{\cc}{{\mathbb{C}}}
\newcommand{\dd}{{\mathbb{D}}}
\newcommand{\ee}{{\mathbb{E}}}
\newcommand{\ff}{{\mathbb{F}}}
\renewcommand{\gg}{{\mathbb{G}}}
\newcommand{\hh}{{\mathbb{H}}}
\newcommand{\ii}{{\mathbb{I}}}
\newcommand{\jj}{{\mathbb{J}}}
\newcommand{\kk}{{\mathbb{K}}}
\renewcommand{\ll}{{\mathbb{L}}}
\newcommand{\mm}{{\mathbb{M}}}
\newcommand{\nn}{{\mathbb{N}}}
\newcommand{\oo}{{\mathbb{O}}}
\newcommand{\pp}{{\mathbb{P}}}
\newcommand{\qq}{{\mathbb{Q}}}
\newcommand{\rr}{{\mathbb{R}}}
\renewcommand{\ss}{{\mathbb{S}}}
\newcommand{\ttt}{{\mathbb{T}}}
\newcommand{\uu}{{\mathbb{U}}}
\newcommand{\vv}{{\mathbb{V}}}
\newcommand{\ww}{{\mathbb{W}}}
\newcommand{\xx}{{\mathbb{X}}}
\newcommand{\yy}{{\mathbb{Y}}}
\newcommand{\zz}{{\mathbb{Z}}}

\newcommand{\AAA}{{\mathbf{A}}}
\newcommand{\BBB}{{\mathbf{B}} }
\newcommand{\CCC}{{\mathbf{C}} }
\newcommand{\DDD}{{\mathbf{D}} }
\newcommand{\EEE}{{\mathbf{E}} }
\newcommand{\FFF}{{\mathbf{F}} }
\newcommand{\GGG}{{\mathbf{G}}}
\newcommand{\HHH}{{\mathbf{H}}}
\newcommand{\III}{{\mathbf{I}}}
\newcommand{\JJJ}{{\mathbf{J}}}
\newcommand{\KKK}{{\mathbf{K}}}
\newcommand{\LLL}{{\mathbf{L}}}
\newcommand{\MMM}{{\mathbf{M}}}
\newcommand{\NNN}{{\mathbf{N}}}
\newcommand{\OOO}{{\mathbf{O}}}
\newcommand{\PPP}{{\mathbf{P}}}
\newcommand{\QQQ}{{\mathbf{Q}}}
\newcommand{\RRR}{{\mathbf{R}}}
\newcommand{\SSS}{{\mathbf{S}}}
\newcommand{\TTT}{{\mathbf{T}}}
\newcommand{\UUU}{{\mathbf{U}}}
\newcommand{\VVV}{{\mathbf{V}}}
\newcommand{\WWW}{{\mathbf{W}}}
\newcommand{\XXX}{{\mathbf{X}}}
\newcommand{\YYY}{{\mathbf{Y}}}
\newcommand{\ZZZ}{{\mathbf{Z}}}

\newcommand{\gA}{{\mathfrak{A}}}
\newcommand{\gB}{{\mathfrak{B}}}
\newcommand{\gC}{{\mathfrak{C}}}
\newcommand{\gD}{{\mathfrak{D}}}
\newcommand{\gE}{{\mathfrak{E}}}
\newcommand{\gF}{{\mathfrak{F}}}
\newcommand{\gG}{{\mathfrak{G}}}
\newcommand{\gH}{{\mathfrak{H}}}
\newcommand{\gI}{{\mathfrak{I}}}
\newcommand{\gJ}{{\mathfrak{J}}}
\newcommand{\gK}{{\mathfrak{K}}}
\newcommand{\gL}{{\mathfrak{L}}}
\newcommand{\gM}{{\mathfrak{M}}}
\newcommand{\gN}{{\mathfrak{N}}}
\newcommand{\gO}{{\mathfrak{O}}}
\newcommand{\gP}{{\mathfrak{P}}}
\newcommand{\gR}{{\mathfrak{R}}}
\newcommand{\gS}{{\mathfrak{S}}}
\newcommand{\gT}{{\mathfrak{T}}}
\newcommand{\gU}{{\mathfrak{U}}}
\newcommand{\gV}{{\mathfrak{V}}}
\newcommand{\gW}{{\mathfrak{W}}}
\newcommand{\gX}{{\mathfrak{X}}}
\newcommand{\gY}{{\mathfrak{Y}}}
\newcommand{\gZ}{{\mathfrak{Z}}}

\newcommand{\gota}{{\mathfrak{a}}}
\newcommand{\gotb}{{\mathfrak{b}}}
\newcommand{\gotc}{{\mathfrak{c}}}
\newcommand{\gotd}{{\mathfrak{d}}}
\newcommand{\gote}{{\mathfrak{e}}}
\newcommand{\gotf}{{\mathfrak{f}}}
\newcommand{\gotg}{{\mathfrak{g}}}
\newcommand{\goth}{{\mathfrak{h}}}
\newcommand{\goti}{{\mathfrak{i}}}
\newcommand{\gotj}{{\mathfrak{j}}}
\newcommand{\gotk}{{\mathfrak{k}}}
\newcommand{\gotl}{{\mathfrak{l}}}
\newcommand{\gotm}{{\mathfrak{m}}}
\newcommand{\gotn}{{\mathfrak{n}}}
\newcommand{\goto}{{\mathfrak{o}}}
\newcommand{\gotp}{{\mathfrak{p}}}
\newcommand{\gotq}{{\mathfrak{q}}}
\newcommand{\gotr}{{\mathfrak{r}}}
\newcommand{\gots}{{\mathfrak{s}}}
\newcommand{\gott}{{\mathfrak{t}}}
\newcommand{\gotu}{{\mathfrak{u}}}
\newcommand{\gotv}{{\mathfrak{v}}}
\newcommand{\gotw}{{\mathfrak{w}}}
\newcommand{\gotx}{{\mathfrak{x}}}
\newcommand{\goty}{{\mathfrak{y}}}
\newcommand{\gotz}{{\mathfrak{z}}}

\newcommand{\krest}{\begin{picture}(14,14)
\put(00,04){\line(1,0){14}}
\put(00,02){\line(1,0){14}}
\put(06,-4){\line(0,1){14}}
\put(08,-4){\line(0,1){14}}
\end{picture}     }

\newcommand{\tret}{{\frac 13}}
\newcommand{\dvet}{{\frac 23}}
\newcommand{\polt}{{\frac 32}}
\newcommand{\polo}{{\frac 12}}

\renewcommand{\leq}{\leqslant}
\renewcommand{\geq}{\geqslant}

\renewcommand{\th}{therefore}
\newcommand{\gr}{gradient}
\newcommand{\Mf}{Morse function}
\newcommand{\iis}{it is sufficient}
\newcommand{\sut}{~such~that~}
\newcommand{\wrt}{~with respect to}
\newcommand{\sm}{\setminus}
\newcommand{\ho}{homomorphism}
\newcommand{\ma}{manifold}
\newcommand{\nei}{neighborhood}

\newcommand{\FR}{{\mathcal{F}}r}
\newcommand{\gt}{{\mathcal{G}}t}

\newcommand{\aand}{\quad\text{and}\quad}
\newcommand{\wwhere}{\quad\text{where}\quad}
\newcommand{\ffor}{\quad\text{for}\quad}
\newcommand{\for}{~\text{for}~}
\newcommand{\iif}{\quad\text{if}\quad}
\newcommand{\iiif}{~\text{if}~}

\newcommand{\Wkr}{W^{\circ}}

\newcommand{\Ker}{\text{\rm Ker }}
\newcommand{\ind}{\text{\rm ind}}
\newcommand{\rk}{\text{\rm rk }}
\renewcommand{\Im}{\text{\rm Im }}
\newcommand{\supp}{\text{\rm supp }}
\newcommand{\Int}{\text{\rm Int }}

\newcommand{\nr}{\Vert}
\newcommand{\smo}{C^{\infty}}

\newcommand{\fpr}[2]{{#1}^{-1}({#2})}
\newcommand{\sdvg}[3]{\widehat{#1}_{[{#2},{#3}]}}
\newcommand{\disc}[3]{B^{({#1})}_{#2}({#3})}
\newcommand{\Disc}[3]{D^{({#1})}_{#2}({#3})}
\newcommand{\desc}[3]{B_{#1}(\leq{#2},{#3})}
\newcommand{\Desc}[3]{D_{#1}(\leq{#2},{#3})}
\newcommand{\komp}[3]{{\bold K}({#1})^{({#2})}({#3})}
\newcommand{\Komp}[3]{\big({\bold K}({#1})\big)^{({#2})}({#3})}

\newcommand{\ran}{\{(A_\lambda , B_\lambda)\}_{\lambda\in\Lambda}}
\newcommand{\rans}{\{(A_\sigma , B_\sigma)\}_{\sigma\in\Sigma}}

\newcommand{\fmin}{F^{-1}}
\newcommand{\vh}{\widehat{(-v)}}

\newcommand{\chart}{\Phi_p:U_p\to B^n(0,r_p)}
\newcommand{\atlas}{\{\Phi_p:U_p\to B^n(0,r_p)\}_{p\in S(f)}}
\newcommand{\flow}{{\VV}=(f,v, \UU)}

\newcommand{\Rn}{\bold R^n}
\newcommand{\Rk}{\bold R^k}

\newcommand{\crr}{p\in S(f)}
\newcommand{\nrv}{\Vert v \Vert}
\newcommand{\nrw}{\Vert w \Vert}
\newcommand{\nru}{\Vert u \Vert}

\newcommand{\obb}{\cup_{p\in S(f)} U_p}
\newcommand{\proob}{\Phi_p^{-1}(B^n(0,}

\newcommand{\stind}[3]{{#1}^{\displaystyle\rightsquigarrow}_
{[{#2},{#3}]}}

\newcommand{\indl}[1]{{\scriptstyle{\text{\rm ind}\leqslant {#1}~}}}
\newcommand{\inde}[1]{{\scriptstyle{\text{\rm ind}      =   {#1}~}}}
\newcommand{\indg}[1]{{\scriptstyle{\text{\rm ind}\geqslant {#1}~}}}

\newcommand{\obbi}{\cup_{p\in S_i(f)}}
\newcommand{\vem}{\text{Vectt}(M)}

\newcommand{\id}{\text{id}}

\newcommand{\xit}{\tilde\xi_t}

\newcommand{\st}[1]{\overset{\rightsquigarrow}{#1}}
\newcommand{\bst}[1]{\overset{\displaystyle\rightsquigarrow}\to{\boldkey{#1}}}

\newcommand{\stexp}[1]{{#1}^{\rightsquigarrow}}
\newcommand{\bstexp}[1]{{#1}^{\displaystyle\rightsquigarrow}}

\newcommand{\bstind}[3]{{\boldkey{#1}}^{\displaystyle\rightsquigarrow}_
{[{#2},{#3}]}}
\newcommand{\bminstind}[3]{\stind{({\boldkey{-}\boldkey{#1}})}{#2}{#3}}

\newcommand{\Tb}{\text{ \rm Tb}}

\newcommand{\VODIN}{V_{1/3}}
\newcommand{\VDVA}{V_{2/3}}
\newcommand{\VM}{V_{1/2}}
\newcommand{\ddd}{\cup_{p\in S_i(F_1)} D_p(u)}
\newcommand{\dddmin}{\cup_{p\in S_i(F_1)} D_p(-u)}
\newcommand{\where}{\quad\text{\rm where}\quad}

\newcommand{\kr}[1]{{#1}^{\circ}}

\newcommand{\vew}{\text{\rm Vect}^1 (W,\bot)}

\newcommand{\Imm}{\text{\rm Im}}
\newcommand{\hrrr}{\text{\rm Vect}^1(M)}
\newcommand{\vemm}{\text{\rm Vect}^1_0(M)}
\newcommand{\ver}{\text{\rm Vect}^1(\RRR^ n)}
\newcommand{\verr}{\text{\rm Vect}^1_0(\RRR^ n)}

\newcommand{\mods}{\vert s(t)\vert}
\newcommand{\exd}{e^{2(D+\alpha)t}}
\newcommand{\exmin}{e^{-2(D+\alpha)t}}

\newcommand{\intt}{[-\theta,\theta]}

\newcommand{\ra}{\rightarrow}
\newcommand{\rau}[1]{\overset{#1}{\ra}}
\newcommand{\la}{\leftarrow}
\newcommand{\lau}[1]{\overset{#1}{\la}}

\newcommand{\Lhat}{\widehat\L}

\newcommand{\wi}{\widetilde}

\newcommand{\wh}{\widehat}

\newcommand{\qs}{\quad\square}

\newcommand{\tens}[1]{\underset{#1}{\otimes}}

\newcommand{\ov}{\overline}

\newcommand{\Hom}{\text{\rm Hom}}

\newcommand{\sphat}{~\widehat{}~}

\newcommand{\pr}{\partial}

\newcommand{\sbs}{\subset}

\newcommand{\Lxi}{\L^-_\xi}

\newcommand{\pa}{\vskip0.1in}

\newcommand{\linia}{\begin{picture}(9,0)
\thicklines
\put(0,0){\line(1,0){490}}
\end{picture}}

\newcounter{c}[section]
\newcounter{d}[section]
\newenvironment{subs}{\vskip0.1in     \noindent
 {\bfseries
\addtocounter{d}{1}\arabic{section}.\arabic{d}. }
\bfseries        }{.\quad }

\newtheorem*{Main}{Main~Theorem}

\renewcommand{\thec}{\thesection.\arabic{c}}

\newcommand{\prop}[1]{\par\refstepcounter{c}%
\vskip0.1in   \noindent {\bf Proposition \arabic{section}.\arabic{c}.}
\hspace{0,2cm}{\it {#1}}\vskip0.1in}

\newcommand{\theo}[1]{\par\refstepcounter{c}%
\vskip0.1in   \noindent {\bf Theorem \arabic{section}.\arabic{c}.}
\hspace{0,2cm}{\it {#1}}\vskip0.1in}

\newcommand{\coro}[1]{\par\refstepcounter{c}%
\vskip0.1in   \noindent {\bf Corollary \arabic{section}.\arabic{c}.}
\hspace{0,2cm}{\it {#1}}\vskip0.1in}

\newcommand{\lemm}[1]{\par\refstepcounter{c}%
\vskip0.1in   \noindent {\bf Lemma \arabic{section}.\arabic{c}.}
\hspace{0,2cm}{\it {#1}}\vskip0.1in}

\newcommand{\defi}[1]{\par\refstepcounter{c}%
\vskip0.1in   \noindent {\bf Definition \arabic{section}.\arabic{c}.}
\hspace{0,2cm}{#1}\vskip0.1in}

\newcommand{\deficit}[2]{\par\refstepcounter{c}%
\vskip0.1in   \noindent {\bf Definition \arabic{section}.\arabic{c} {#1}.}
\hspace{0,2cm}{ {#2}}\vskip0.1in}

\newcommand{\cons}[1]{\par\refstepcounter{c}%
\vskip0.1in   \noindent {\bf Construction \arabic{section}.\arabic{c}.}
 \quad{\it {#1}}\vskip0.1in}

\newcommand{\theocit}[2]{\par\refstepcounter{c}%
\vskip0.1in   \noindent {\bf Theorem \arabic{section}.\arabic{c} {#1}.}
\hspace{0,2cm}{\it {#2}}\vskip0.1in}

\newcommand{\rema}[1]{\par\refstepcounter{c}%
\vskip0.1in   \noindent {\bf Remark \arabic{section}.\arabic{c}.}
\hspace{0,2cm}{#1}\vskip0.1in}

\title[Asymptotics of Morse numbers]
{On the asymptotics of Morse numbers
of finite          covers of   manifolds       }
\author{A.V.Pajitnov}
\address{Universit\'e de Nantes, Facult\'e des Sciences,
2, rue de la Houssini\`ere, 44072, Nantes Cedex}
\email{pajitnov@math.univ-nantes.fr}
\begin{abstract}

Let $M$ be a closed connected manifold.
We  denote by $\MM(M)$ the Morse number of $M$,
that is, the minimal possible number of critical points of
a Morse function $f$ on $M$.
M.Gromov posed the following question:
 Let $N_k, k\in \NNN$  be a
sequence of manifolds, such that each $N_k$ is an $a_k$-fold cover
of  $M$ where $a_k\to\infty$
as $k\to\infty$.
What are the asymptotic properties of the sequence
$\MM(N_k)$ as $k\to\infty$?

In this paper we study the case $\pi_1(M)\approx\ZZZ^m, \dim M\geq 6 $.
Let               $\xi\in H^1(M,\ZZZ), \xi\not=0  $.
Let $M(\xi)$ be the infinite cyclic cover
corresponding to $\xi$, with generating covering
translation $t:M(\xi)\to M(\xi)$.
Let
   $M(\xi,k)$
be the
quotient
$M(\xi)/t^k$.
We prove that
 $\lim_{k\to\infty} \MM(M(\xi,k))/k$ exists.
For $\xi$ outside  a subset $\gM\subset H^1(M)$ which is the
union of a finite family of hyperplanes,
 we obtain the asymptotics of
$\MM(M(\xi,k))$ as $k\to\infty$ in terms of
homotopy invariants of $M$ related to the
 Novikov
homology of $M$.
It turns out that the limit above does not depend on $\xi$
(if $\xi\notin \gM$).
Similar results hold for the stable Morse numbers.
Generalizations for the case of non-cyclic coverings are obtained.

\end{abstract}
\maketitle
\section*{Introduction and the statement of the result}

Let $M$ be a closed connected smooth     manifold. Denote
by $\MM(M)$ the Morse number of $M$, that is, the
minimal possible number of critical points of a Morse function on
$M$. In the case $\pi_1(M)=0, \dim M\geq 6$, this number
is easily computable  in terms of homology of $M$ (see   \cite{Sm}).
In the case of arbitrary fundamental group (even for $\dim M\geq 6$),
the number $\MM(M)$ is very difficult to compute:
it depends on the simple homotopy
type of $M$, the relevant algebraic constructions
are rather complicated, and it is not easy to extract the
needed numerical invariant (see \cite{Sh1}, or  \cite{Sh2}, Ch. 7).

\pa
M.Gromov posed the following question:
\vskip0.2in
{\it Let $N_k, k\in \NNN$  be a
sequence of manifolds, such that each $N_k$ is an $a_k$-fold cover
of the manifold $M$ where $a_k\to\infty$
as $k\to\infty$.
What are the asymptotic properties of the sequence
$\MM(N_k)$ as $k\to\infty$?
\vskip0.1in
}
In the present article we study the problem for
$\pi_1(M)$
free abelian and
$\dim M\geq 6$.
To formulate our results, we need some terminology from algebra.
Denote $\ZZZ[\ZZZ^m]$ by $\L$. Let
$$C_*=\{0\la C_0\lau{\partial_1} C_1 . . . \lau{\partial_k} C_k\la 0\}$$
be a free finitely generated $\L$-complex. Denote by $B_i(C_*)$ the rank of the module
$H_i(C_*)\underset{\L}{\otimes}\{\L\}$
over the field of fractions $\{\L\}$.
Denote by $B(C_*)$ the sum of all
$B_i(C_*)$.
 Consider  now the homomorphism
$\pr_{i+1}:C_{i+1}\to C_i$, and let $d=\rk C_i$. Recall that
the Fitting invariant $\FF_t$ of the homomorphism
$\pr_{i+1}$
(see e.g. [4], p.278) is the ideal of $\L$ generated by
the
$(d-t)\times(d-t)$ subdeterminants of the matrix of
$\pr_{i+1}$
(for $t\geq d$ one sets
$\FF_t=\L$ by definition).
We shall denote the sequence
$\FF_0\sbs\FF_1\sbs .. \sbs \FF_d$
of the Fitting invariants by $F(i)$.

Define the {\it reduced Fitting sequence }
for $\partial_{i+1} $ to be the sequence
\begin{equation*}
 \FF_s\subset . . . \subset \FF_r\tag*{$FR(i)$}
\end{equation*}
where $\FF_s$, respectively $\FF_r$, is the first, respectively the last, term
of the Fitting sequence $F(i)$, not equal to $0$, respectively to $\L$.
The sequence $F(i)$ is not a homotopy invariant
of $C_*$, but the sequence $FR(i)$ is (see
e.g. \cite{Sh2}, Ch.4,\S 2).
 We say that an ideal $J$ of $\L$ is
{\it numerically prime} if there is no
number $l\in\ZZZ, l\not=\pm 1$, \sut~ every $R\in J$ is divisible by
$l$, and we denote by
$Q_i(C_*)$ the number of ideals in the sequence
$FR(i)$ which are {\it not} numerically prime.
Denote by $Q(C_*)$ the sum of all $Q_i(C_*)$.

A subgroup $G\subset\ZZZ^m$ will be called an
{\it integral hyperplane} if it is a direct summand of $\ZZZ^m$
of rank $m-1$.

Now let $M$ be a closed connected manifold,
 $\pi_1(M)\approx\ZZZ^m, m\geq 1$.
It is convenient to set $m=n+1,n\geq 0$.
 For every
non-zero $\xi\in H^1(M)$ there is a unique connected infinite cyclic
covering $\PP_\xi  :M(\xi)\to M$ \sut
$\PP_\xi ^*(\xi)=0$. Denote by
$M(\xi,k)\to M$ the $k$-fold
cyclic covering of $M$
 obtained from $\PP_\xi   $.
Let $C_*(\widetilde M)$ be the
cellular chain complex of the universal cover
  $\wi        M$.
  We shall abbreviate
$B(C_*(\wi M))$ to $B(M)$ and $Q(C_*(\wi M) )$ to $Q(M)$.

\begin{Main}
Let $\dim M\geq 6, \pi_1(M)\approx\ZZZ^{n+1}, n\geq 0  $. Then:
\begin{enumerate}\item
For any non-zero
$\xi\in H^1(M)$
the limit
$\lim_{k\to\infty} \MM(M(\xi,k))/k$
exists.

\item  There is a subset $\gM\subset H^1(M)$
which is  a finite union of integral hyperplanes in $H^1(M)$,
and for every non-zero $\xi\notin \gM$ there is
a real number $a$ \sut~ for every $k\in\NNN$ we have
$$
k(B(M)+2Q(M))-a\leq\MM(M(\xi,k))\leq k(B(M)+2Q(M))+a
$$
\end{enumerate}       \end{Main}
{\it Remarks:\quad}
1) A similar result holds for the stable Morse numbers of $M$, see \S 5.

2) The limit $\lim_{k\to\infty} \MM(M(\xi,k))/k  $ will be denoted by
$\mu(M,\xi)$. The second point of the Main Theorem implies that
for a "generic" cohomology class $\xi$
we have $\mu(M,\xi)=B(M)+2Q(M)$.

3) Denote by $\MM_i(M)$ the minimal number of critical
points of index $i$ of a Morse function on $M$. The methods of
the present paper allow also to prove that (under the assumptions of
the Main Theorem) the limit
$\lim_{k\to\infty}\MM_i(M(\xi,k))/k$ exists, and that for all $\xi$,
except those belonging to a finite union of integral hyperplanes,
there is a real number $a$ \sut~ for every natural $k$ we have
$$k(B_i(M)+Q_i(M)+Q_{i-1}(M))-a
\leq \MM_i(M(\xi,k)) \leq
k(B_i(M)+Q_i(M)+Q_{i-1}(M))+a.
$$

4) The numbers $B_s(M), Q_s(M)$ are closely related to
the Novikov homology of $M$. Namely, $B_s(M)$ is equal to the
Novikov Betti number $b_s(M,\xi)$ \cite{Nov1}
for every non-zero class $\xi\in H^1(M)$ (note that $B_s(M)$ is also equal
to the $L^2$-Betti number $b_s^{(2)}(M)$).
Further, for every non-zero $\xi\notin \gM$
we have
$Q_s(M)\leq q_s(M,\xi)$ where
$q_s(M,\xi)$ is the Novikov torsion number  \cite{Nov1}( that follows from
Remark 2.6 and Proposition 3.3 of the
present paper).

\vskip0.1in
The proof  is outlined as follows. Assume that $\xi\in H^1(M)$ is
indivisible. Let $f:M\to S^1$ be a Morse map, representing
$\xi$, and let $V=f^{-1}(\l)$ be a regular level surface of $f$.
We can assume that $V$ is connected and that $\pi_1(V)\to\pi_1(M)$
is an isomorphism onto $\Ker\xi$.
Cut $M$ along $V$, and obtain  a cobordism $W$, \sut~
the boundary $\pr W$ has two connected components $\pr_0W$ and $\pr_1 W$,
each diffeomorphic to $V$.
The cyclic cover  $M(\xi)$ is the union of a countable family of copies of
$W$ glued successively
to each other. The union $W_k$  of $k$ successive copies
is a cobordism. Its boundary $\pr W_k$
has two connected components
$\pr_0W_k$ and $\pr_1 W_k$,
each diffeomorphic to $V$
(see \S 4 for details).
We show that $\MM(M(\xi,k))$ and $\MM(W_k,\pr_0 W_k)$
have the same asymptotics
  as $k\to\infty$
(see \S 4).
Further, $\MM(W_k,\pr_0 W_k)$ is equal to the Morse number
of the $\ZZZ[\ZZZ^n]$-complex $C_*(\wi W_k, \wi{ \pr_0 W_k})$,
see \S 1
 for definitions.
It turns out that the asymptotic behaviour
of this Morse number (as $k\to\infty$)
depends only on
 the chain homotopy type of
$C_*(\wi{M(\xi)})$ (moreover it depends only on the
Novikov completion of this complex).
 The definition and the properties of the corresponding
invariant of chain complexes
are the subject of \S\S 1 -  3
of the paper. These sections  are purely algebraic.
It follows from the author's earlier result \cite{PaLor}, that for
$\xi$ outside  a finite union of integral hyperplanes
in $H^1(M)$, the Novikov-completed chain homotopy type
of $C_*(\wi{M(\xi)})$ is easily computable. (This is the subject
of the second half of \S 2 and of \S 3.)
This leads to the effective computation of the
asymptotics presented in the main theorem.

                 \vskip0.1in

I am grateful to M.Gromov for a  stimulating discussion
on the subject. He suggested in particular, that asymptotically
the numbers $\MM(N)$ above should be related
to Novikov numbers. He indicated also that the Morse number
$\MM(M(\xi,k))$
should have the same asymptotics as the
 Morse number of the            pair         $(W_k,\pr_0 W_k)$
(see Prop. \ref{varb}  of the present paper).

\section{Morse numbers of chain complexes}

In this section we define the notion of the Morse number
for arbitrary chain complexes over
 $\ZZZ[\ZZZ^n]$ and we develop some
basic
 properties
of these numbers.   We assume  that the reader is familiar with
\S 3 of   \cite{Sh1}  and with \S 1 of \cite{PaSur}.
We denote $\ZZZ[\ZZZ^n]$ by $R$.

{\it Terminological remark.\quad} Let $A_*, B_*$ be chain complexes.
We shall denote the
 chain maps from $A_*$ to $B_*$ as follows: $f_*:A_*\to B_*$,
so that $f_k$ is a homomorphism $A_k\to B_k$.

\defi{
                   An  $R$-{\it complex} is a chain complex
$\{0\la C_0\la C_1 . . . \la C_k \la 0  \}$
of finitely generated $R$-modules.
The {\it length }  $l(C_*)$ of an $R$-complex $C_*$
is the maximal number $l$ \sut $C_l\not=0$.
An $R$-complex $  C_*$ is called a {\it free $R$-complex   }
(or simply     {\it   $f$-complex}     )
if every $C_i$ is a free finitely generated module over $R$.
         }

\deficit{\cite{Sh1}}{
Let $C_*$ be an $f$-complex over  $R$.
The minimal possible number of free generators of an $f$-complex
$ D_*  $, having the same homotopy type as $ C_* $,
is called the {\it Morse number of} $C_*$ and denoted
by $\MM(C_*)$
(or by $\MM_R(C_*)$, if we want to stress
the base ring).
}

One of the consequences of the Quillen-Suslin    theorem (\cite{Quil},\cite{Sus})
is that $R$ is an $s$-ring, that is, every projective $R$-module
is free (see \cite{Lam}, Ch.5, \S 4). $R$ is also an $IBN$-ring, that is,
the number of free generators of a free module is uniquely determined.
Therefore,
 in the homotopy type of every $f$-complex over $R$
there exists a minimal chain complex, that
is, a complex $D_*$ such that the number of
free generators of $D_*$ in each dimension is minimal
over all the free complexes in this homotopy type
(see \cite{Sh1}, Th.3.7).

\defi{
Let $A_*$ be an $R$-complex.
We call a {\it free model} of 
$A_*$ a free $R$-complex $A'_*$ together with a
chain map
 $\a_*: A'_*\to A_*$ which is epimorphic and
 induces an isomorphism in homology.
                     \footnote{Sometimes we shall say  (by abuse of
terminology) that the complex $A'_*$ itself   is a free model of $A_*$.}
If
 $\a_*: A'_*\to A_*$,
 $\b_*: B'_*\to B_*$
are free models, 
and $f_*:A_*\to B_*$
is a chain map,
then
a chain map  $F_*:A'_*\to B'_*$
is called {\it covering} of $f$
if    $\b_* F_* = f_*\a_*$. Similar terminology
is accepted for chain homotopies.
         }

\lemm{
Let $A_*$ be an $R$-complex. Then
there is a free model of $A_*$ and
\begin{enumerate}
\item Every chain map $A_*\to B_*$ admits a covering
\wrt~ any free models of $A_*$ and $B_*$.
\item Let $h_*:A_*\to B_{*+1}$ be a chain homotopy
from $f_*$ to $g_*$, and $F_*, G_*$ be coverings of
$f_*$, respectively $g_*$,
\wrt~ some free models of $A_*, B_*$.
 Then there is a chain homotopy
$H_*$ from $F_*$ to $G_*$, covering $h_*$.
\item Two free models of a complex $A_*$ are homotopy equivalent.
\end{enumerate}
         }

{\it Proof.\quad} To prove the existence of a free model,
we proceed
by induction in the length of $A_*$.
If $l(A_*)=0$, then it follows from the fact that  every
finitely generated module over $R$ has a free finite resolution of finite length.
To make the induction step, it suffices to construct a free model
for a complex  of the type
$C_*=
\{  0\la A_0\lau{\partial_1} C_1\lau{\partial_2} C_2 . . . \lau{\partial_n}
C_n\la 0\}$,
where $C_i$ are free finitely generated
 modules and $A_0$ is a finitely generated module.
Let $B_*=       \{0\la       A_0
\lau{\e}
E_0
\lau{d_1}
E_1
\lau{d_2} . . .\}            $
be  a finite free resolution of $A_0$.
There is a
chain
map
$\phi_*:C_*\to B_*$, \sut $\phi_0=\text{\rm id}$.
Define now  an $R$-complex
$$F_*=\{0\la E_0
\lau{D_1}
C_1\oplus     E_1
\lau{D_2}
C_2\oplus E_2
. . .        \}$$
setting
$D_1(c_1,e_1)=\phi_1(c_1)+d_1(e_1)$ and
 $D_i(c_i,e_i)= (\partial_i(c_i), d_i(e_i)+
(-1)^{i+1}\phi_i(c_i))$ for $i\geq 2$.

Define further a map $\g_*  :F_*\to C_*    $
to be the projection $(x,y)\mapsto x$      when $*\geq 1$ and set
$\g_0=\e$.
It is easy to check that
$F_*$ is indeed an
$f$-complex, and that $\g_*$ is a free model.
The points (1) and  (2) of our lemma
are proved by a standard homological algebra argument;
(3) follows from (2).
$\quad \square$

\defi{
 The {\it Morse number $\MM(C_*)$} of a complex $C_*$
is the Morse number of any of its
free models.    
         }

\prop{
Let
$0\la A_*\la B_*\la C_*\la 0$
be an exact sequence of $R$-complexes. Then
1)   $\MM(B_*)\leq \MM(A_*)+\MM(C_*)$, and
2)   $\MM(A_*)\leq \MM(C_*)+\MM(B_*)$.
         }
\label{pryamsum}

{\it Proof.\quad}
1)   The following lemma reduces the assertion to the case of
free $R$-complexes.
\lemm{
Let
$0\la A_*\la B_*\la C_*\la 0    $
be an exact sequence of $R$-complexes.
Then there is a commutative diagram
$$
\begin{CD}
0  @<<<  A'_*  @<<<    B'_*        @<<<      C'_*   @<<< 0  \\
@VVV    @V{\a_*}VV    @V{\b_*}VV            @V{\g_*}VV      @VVV        \\
0  @<<<  A_*  @<<<    B_*        @<<<      C_*   @<<< 0
\end{CD}
$$
where $\a_*,\b_*,\g_*$ are free models.
         }
{\it Proof of the lemma.~}
Let $g'_*:   C'_*\to B'_*$ be a covering of $C_*\to B_*$
\wrt~ some free        models $B'_*, C'_*$.
We can  assume that $g'_*$ is a monomorphism onto
a direct summand (the proof repeats almost verbally the proof of
Lemma 1.8 from \cite{PaSur} and will be omitted).
 Now, setting $A'_*=B'_*/\Im C'_*$, we obtain the
first line of the commutative diagram above.
$     \square$

For   the case of free complexes the assertion
follows from the next one.
\lemm{
Let
$0  \la  A_*  \la    B_*        \la      C_*   \la 0$
be an exact sequence of free $R$-complexes.
Then there is an exact sequence
$0  \la  A'_*  \la    B'_*        \la      C'_*   \la 0$
of free $R$-complexes \sut~
$A'_*\sim A_*,~B'_*\sim B_*,~C'_*\sim C_*$, and $A'_*, C_*'$ are minimal.
         }
The proof of this lemma is an exercise in
 the    theory of
minimal complexes
(\cite{Sh2},\S 4),
 and will be left to  the reader.
$\quad\square$

To prove 2) observe that there is an     exact sequence
$\ZZ=        \{   0  \la  \Sigma C_*  \la
   D_*        \la      B_*   \la 0      \} $
where $D_*$ is the mapping cone of $j_*$, and $\Sigma C_*$
is the suspension of $C_*$.
Now apply the point 1) to the sequence $\ZZ$. $\qs$

In some cases the first inequality of the preceding proposition
turns to    equality. We shall say that a complex $C_*$
{\it is concentrated in dimensions $[k,r]$} if
$C_i=0$   for $i<k$ and for $i>r$.
We denote by $F(i,s)_*$ the chain complex $\{0\la R^s\la 0\}   $ 
concentrated
in dimensions $[i, i]$.
\lemm{
1) For every f-complex $C_*$ we have
$\MM(C_*\oplus F(i,s)_*  )=\MM(C_*)+s$.
2) Let $C_*, D_*$ be f-complexes, concentrated respectively in
dimensions $[a,b]$, and $[b,c]$. Then $\MM(C_*\oplus D_*)=
\MM(C_*)+\MM(D_*)$.
}
\label{razdel}
The proof of this lemma is easily obtained from  V.V.Sharko's
criterion of minimality of chain complexes (see \cite{Sh1}, Lemma 3.6).
$\qs$

\section{A
numerical
 invariant of free chain $R((t))$-complexes}

We 
denote 
       $ \ZZZ[\ZZZ^n]$ by $R$
(as in  the previous section).
Let us start with a free $R[[t]]$-complex $A_*$.
 For  $k\in\NNN$
denote by $A[k]_*$ the     free $R$-complex $A_*/t^kA_*$,
and denote its Morse number
by $\m_k(A_*)$,
so  
$\mu_k(A_*)=\MM_R (A[k]_*)$.
Note that $\m_k(A_*)+\mu_l(A_*)\geq \m_{k+l}(A_*)$.
(Indeed, consider the short exact sequence
$0\la A[k]_* \la A[k+l]_* \la A[l]_* \la 0$
and apply Proposition 1.6.)
Therefore the sequence $\{\m_k/k\}_{k\in\NNN}$ has a limit
(see \cite{PoSe}, ex. 98) which will be denoted by $\s(A_*)$.
It is clear that
 $\s(A_*)$ is
a chain
  homotopy invariant of  $A_*$.

Now we shall consider free complexes over the ring
$R((t))=\s^{-1}R[[t]]$ where $\s$ is the multiplicative set
$\{t^l\mid l\in\NNN\}$.
Let $C_*$ be such a complex. We say that a chain subcomplex
$D_*\subset C_*$ is a {\it basic subcomplex}
if                                           \quad
1)    $    D_*$ is a free $R[[t]]$-complex, \quad and
2)  $\s^{-1}D_*=C_*$.
It is clear that each free complex $C_*$ over $R((t))$
  has  basic subcomplexes.

\prop{
Let $C_*$ be a
free
$R((t))$-complex. Then the number
$\s(D_*)$ is the same for every basic subcomplex
$D_*\subset C_*$.
         }
\label{chislosigma}
{\it Proof.\quad}
Let $D_*, F_*$ be
basic subcomplexes.
The Noetherian property of $R[[t]]$ and the condition (2)
in the definition of a basic subcomplex imply immediately that
there is $s\in\NNN$ \sut~
$t^sF_*\subset D_*$. Since $\s(D_*)=\s(t^kD_*)  $
we can assume that $t^sF_*\subset D_*\subset F_*$.
Now for every $l\in\NNN$
we obtain two exact sequences of finitely generated chain 
complexes over $R$.
\begin{gather}
0\la F_*/D_* \la F_*/t^lD_* \la D_*/t^lD_* \la 0  \label{one}\\
0\la F_*/t^lD_* \la F_*/t^{l+s}F_* \la t^lD_*/ t^{l+s}F_* \la 0  
 \label{two}
\end{gather}

Applying Prop. 1.6 we deduce from  (\ref{one}) and (\ref{two})
that $\m_{l+s}(F_*)\leq C+\m_l(D_*)$ where $C$ does not depend on $l$.
This implies easily that $\s(F_*)\leq \s(D_*)$; by symmetry
we obtain $\s(F_*)=\s(D_*)$. $\quad\square$

Now we can define an invariant of $R((t))$-complexes.
Namely, if $C_*$ is a free $R((t))$-complex, we set
$s(C_*)=\s(D_*)$ where $D_*$ is any basic subcomplex
of $C_*$.
The number $s(C_*)$ depends only on the homotopy type
of the $R((t))$-complex $C_*$. Indeed, a version of the
Cockroft-Swan theorem (\cite{PaSur}, Prop. 1.7)   shows that it is sufficient
to check that $\s(C_*)$ does not change when we add
to $C_*$ a complex of the form
$\{  0\la R((t))\lau{\text{\rm id}} R((t))\la 0\}   $. But this is obvious.

For some free $R((t))$-complexes the asymptotic properties of the 
Morse numbers
are still better.
We shall say that a sequence $a_k$ of real numbers is
{\it asymptotically linear}
if $\exists C,\a,\forall k : \a k-C\leq a_k\leq \a k+C$.
We shall say that a free $R((t))$-complex $C_*$ is of
{\it asymptotically linear growth}
(abbreviation: $aslg$) if for some
basic subcomplex $D_*\subset C_*$
the sequence $\mu_k(D_*)$ is asymptotically linear.
Similarly to the proof of Proposition \ref{chislosigma},
 one can show that in an $aslg$-complex
{\it  every} basic subcomplex $D'_*$ has an asymptotically linear
sequence $\mu_k(D'_*)$.
Note also   that the property of being $aslg$ is homotopy invariant.
We do not know if every $R((t))$-complex is
$aslg$, but  we shall prove that
every complex of a certain class appearing in our geometrical
setting            is $aslg$, and we shall calculate its $s$-invariant.
 We need some definitions.
{\it A monomial } of $R$ is an element of the form
$ag$ where
$a\in\ZZZ$, and
 $g\in\ZZZ^n$.  Let
$Z=z_kt^k+...+z_lt^l\in R[t,t^{-1}]$ where $l,k\in\ZZZ,
k\leq l$, and  $z_k, z_l\not=0   $.
We shall say that                                    $Z$ is:
\begin{itemize}
\item {\it monic} if $z_k=\pm g, g\in\ZZZ^n $ \quad (
Our terminology differs here from the standard one.)
\item    {\it numerically prime}
if it is not divisible by an integer not equal to $\pm 1$.
\item     {\it    special     } if
each $z_i$ is a monomial in $R$.
\end{itemize}

 We denote
$R((t))$ by $\LL$.

\defi{
Let $C_*$ be a complex over $ \LL$. We shall say
that $C_*$ is {\it of principal type}
if
for every $i$  an isomorphism
\begin{equation}                 \label{princ}
H_i(C_*) \approx
 \big(
\bigoplus\limits_{j=1}^{b_i}\LL          \big)
\oplus
\big(
\bigoplus\limits_{s=1}^{q_i} \LL/a_s^{(i)}\LL
\big)
\end{equation}
 is fixed,
and for every $i,s$:
 1)      $a_s^{(i)}\in R[t, t^{-1}]$ and $a_s^{(i)}$ is    special,
non-zero and not monic
  2)
 $a_s^{(i)}\mid a_{s+1}^{(i)}$.
         }
\label{comspec}

For a complex $C_*$ of  principal type
we denote by $\k_i$ the number of those     polynomials
$a_s^{(i)}$ which are not numerically prime.

\theo{
Let $C_*$ be a free $\LL$-complex of principal type.
Then $C_*$ is of asymptotically linear growth, and
 $s(C_*)=\sum_ib_i+2\sum_i\k_i$.
}
\label{comprinc}
{\it Proof.\quad}
We can assume
that all the elements $a^{(i)}_s$ in   (\ref{princ})
are of the form $z_0+...+z_kt^k$ where $z_0\in\ZZZ, z_0\not=0$.
Denote by $\FF(i)_*$ the free $\LL$-complex $\{0\la \LL^{b_i}
\la 0 \}$ concentrated in     dimensions $[i,i]$.
For $\rho\in\LL$ and $i\in\NNN$, denote
by $\tau(\rho,i)_*$   the free complex
$\{0\la \LL \lau{\rho} \LL \la 0\}$ concentrated
in  dimensions
$[i,i+1]$.
Note that
if $\rho\in R[[t]]$ then
 $\tau(\rho,i)_*$ has a standard basic subcomplex
$\{0\la R[[t]] \lau{\rho} R[[t]]\la 0\}$ which will be denoted
by $\tau'(\rho,i)_*$.

For a given $i$
denote by 
$\pi$ (resp. by $\nu$) the set of all $s$ \sut~ $a_s^{(i)}$
is numerically prime (resp. {\it not } numerically prime). 
Set
$$\TT\PP(i)_*=\bigoplus_{s\in\pi} \tau(a_s^{(i)}, i)_*, \quad 
\TT\NN(i)_*=\bigoplus_{s\in\nu} \tau(a_s^{(i)}, i)_*, \quad 
\TT(i)_*=\TT\PP(i)_*\oplus\TT\NN(i)_* .$$

(Morally, $\FF(i)_*$ corresponds to the free part of
the homology $H_i(C_*)$, and $\TT(i)_*$ to the torsion part.)
The complexes
$\TT(i)_*, \TT\PP(i)_*, \TT\NN(i)_*$ have basic subcomplexes
\break
$\TT'(i)_*, \TT\PP'(i)_*, \TT\NN'(i)_*$
which are obtained as direct sums of the     corresponding
complexes
$\tau'(\rho,i)_*$.

Lemma 5.1  of \cite{PaFree}
implies that $C_*$ is homotopy equivalent to
the direct sum (over all $i$)     of the complexes
$\FF(i)_*\oplus \TT(i)_*$.
We call this direct sum {\it principal model} for $C_*$.
Applying
successively
Lemma 1.9,
it is easy to deduce our theorem     from
the next lemma.
\lemm{
For every $i,k$ we have:
(1)   $     \mu_k(\TT\NN'(i)_*)=2k\k_i$, \quad
(2)   $      \mu_k(\TT'(i)_*)\geq 2k\k_i$. \quad
(3) For every $i$ the sequence $\{      \mu_k(\TT\PP'(i)_*)\}_{k\in\NNN}$
 is bounded.
         }
\label{chislas}
{\it Proof.\quad }
1) Fix some $i$.
The condition 2) from the definition \ref{comspec}
implies that there is a prime number $p$ \sut~
every  polynomial $a_s^{(i)}$
which is not numerically prime
 is divisible  by $p$.
Abbreviate
 $\TT\NN'(i)_*$  to
$L_*$; the inequality $\mu_k(L_*)\leq 2k\k_i$ is immediate.
 To  prove the inverse inequality consider
an
 $\FFF_p$-complex
$L[k]_*\tens{R} \FFF_p$
(where $\FFF_p$ is considered as $R=\ZZZ[\ZZZ^n]$-module via
the trivial $\FFF_p$-representation of  $\ZZZ^n$).
It is obvious that $\MM(L[k]_*)$ is not less than
$\dim H_*(L[k]_*\otimes\FFF_p)$ which equals $2k\k_i$. A 
similar argument proves
the point 2).
To prove 3) it suffices to show that
if $\rho=a_0+a_1t+...+a_rt^r\in R[t]$ is special, numerically prime,
and has $a_0\not=0$, then $\mu_k(\tau'(\rho,i)_*)$ is bounded.
Write
$a_j=A_jg_j$ with $ g_j\in\ZZZ^n$ and $ A_j, a_0\in\ZZZ$, and note
 that $a_0, A_1,...,A_r$ are relatively prime.
Abbreviate
$\tau'(\rho,i)_*$ to $S_*$.   If $k\geq 1$
then
 $S[k+r ]_*$
is a free $R$-complex
of the form
$\{0\la R^{k+r} \lau{\AA_k} R^{k+r} \la 0\}$
where $\AA_k$
is
the following
  $(k+r) \times(k+r) $-matrix
\begin{equation}
\AA_k=
\begin{pmatrix}
a_0    &  0   &   \dots  & 0 & 0       \\
a_1    &  a_0 &   \dots  & 0 & 0       \\
\vdots  & \vdots &   \ddots  & \vdots & \vdots   \\
a_r     & a_{r-1} &  \dots  &    0 & 0                \\
0       & a_r     &  \dots   & a_0 & 0            \\
\dots   & \dots   &\dots    & a_1 & a_0
\end{pmatrix}
\end{equation}
Denote by $\BB_k$ the matrix formed by the first $k$ columns
of the matrix $\AA_k$, and denote by
$I_k$ the ideal of $R$ generated by all $k\times k$ subdeterminants
of $\BB_k$.
 The point 2) of the
next
 lemma
implies  our assertion.

\lemm{
1)    $I_k=R$ ;\quad 2)
The submodule $\AA_k(R^{k+r})$ of $R^{k+r}$ contains
a direct summand of $R^{k+r}$, which is a free module of
rank $k$.
}
\label{split}
{\it Proof.\quad}
 The point 2) follows from the point 1)
 by a standard
argument based on the Quillen-Suslin       theorem (we leave the details
to the reader).
Proceeding  to the proof of the point 1), note
 that $a_0^k\in I_k$.
Therefore we can assume that $a_0\not= \pm 1$.
 It suffices to show that for every
prime number $p$ from the prime decomposition of $a_0$
there is an element $C\in R$ \sut~ $1+pC\in I_k$.
To show this, recall that the numbers
$A_j$ are relatively prime. So there is
$i$ \sut~ $p\mid A_j$ for $j<i$, and $p\nmid A_i$.
Consider the
$k\times k$
subdeterminant of the matrix $\BB_k$ formed by
all the columns and by the lines from $i+1$ to $i+k$.
The terms of the principal diagonal are all equal to
$A_i$, and the terms above the diagonal
are divisible by $p$. Therefore this subdeterminant  equals
$Q+pC$ where $Q$ is a monomial of the form
$Q=qg$ with $(q,p)=1   $,
and  Lemma 2.5 follows. This finishes the proof of Theorem 2.3.
 $\quad\square$

\rema{
Let $L=\ZZZ[t, t^{-1}], \hat L =\ZZZ((t))$. The homomorphism
$\e:\ZZZ^n\to \{ 1\}$ extends to   ring homomorphisms
$e:R[t, t^{-1}]\to L$ and $\hat e :\LL\to\hat L$.
Therefore for every $\LL$-complex $C_*$ we can form an $\hat L$-complex 
$\overline{C_*}=C_*\tens{\LL}\hat L$.
 Assume that $C_*$ is of principal type.
Using the homotopy equivalence
$C_*\sim \oplus_i (\FF(i)_*\oplus \TT(i)_*)$ from the proof of Theorem
2.3, it is easy to see that $\overline{C_*}$ is also of principal type
and
\begin{equation}
H_i(\overline{  C_*}) \approx
 \big(
\bigoplus\limits_{j=1}^{b_i}\hat L      \big)
\oplus
\big(
\bigoplus\limits_{s=1}^{q_i} \hat L /\a_s^{(i)}\hat L
\big)
\end{equation}
with $\a_s^{(i)}=e(a_s^{(i)})$. Since
$a_s^{(i)}$ are special and not monic, $b_i$ and $q_i$ are equal
respectively to the rank and to the torsion number of $H_i(\overline{C_*})$
over the principal ring $\hat L$.
It is easy to see that the above decomposition satisfies
Definition 2.2, therefore $\overline{C_*}$ is $aslg$. Further,
$a_s^{(i)}$ is numerically prime if and only if $\a_s^{(i)}$
is, and this  implies
$s(C_*)=s(\overline{C_*})$. }\label{Rema}

\section{A numerical invariant $S(C_*,\xi)  $}

In this section  $\L=\ZZZ[\ZZZ^{n+1}]$,
 $C_*$ is a free
$\L$-complex,
 and $\xi:\ZZZ^{n+1} \to\ZZZ$ is a non-zero  homomorphism.
We    define a numerical invariant
$S(C_*,\xi)$.
For the cohomology classes $\xi$ outside  a finite union
of integer hyperplanes we calculate $S(C_*,\xi)$
in terms of the reduced Fitting sequences of the boundary operators 
of $C_*$ (the mentioned finite union of integer hyperplanes
depends on $C_*$).
 An element $z\in\L$
 is called {\it $\xi$-monic} if
$z=\pm g+z_0$ where $g\in\ZZZ^{n+1}$ and $\supp z_0
\subset\{h\in\ZZZ^{n+1}\mid\xi(h)<\xi(g)\}$.
An element $z$ is called {\it $\xi$-special}
if any two different elements $a,b\in\supp z$ satisfy $\xi(a)\not=\xi(b)$.
We denote by $S_\xi$ the multiplicative subset of all $\xi$-monic
polynomials, and we denote by $\L_{(\xi)}$ the localization $S^{-1}_\xi\L$.

\defi{
A subset $X\subset\ZZZ^k$ will be called
{\it small} if it is a finite union of integer
hyperplanes.
         }
\theocit{(\cite{PaLor}, Th. 0.1)}{
There is a small
subset $\gN\subset \Hom(\ZZZ^{n+1},\ZZZ)$
\sut~ for every
$\xi\notin \gN$ and every $p$ we have:
\begin{equation} \label{Lor}
S_\xi^{-1} H_p(C_*) \approx
\bigg(\bigoplus\limits_{i=1}^{b_p(C_*,\xi  )}
\L_{(\xi)}
\bigg)
\oplus
\bigg(\bigoplus\limits_{j=1}^{q_p(C_*,\xi    )}
\L_{(\xi)}/a_j^{(p)}\L_{(\xi)}
\bigg)
\end{equation}
where  $a_j^{(p)}\in\L$
are non-zero and not $\xi$-monic elements of $\L$ (depending on $\xi$),
and
 $a_j^{(p)}\mid a_{j+1}^{(p)}$. $\qs$
         }

\pa
{\bf Sketch of the proof of Theorem 3.2.}
We shall recall here the basic idea of the proof of 3.2 
following \cite{PaFree}
and
 \cite{PaLor}, see  \cite{PaLor} for the full proof.
Let 
$\xi:\ZZZ^{n+1}\to\RRR$ be a non-zero homomorphism.
Similarly to the above, we define the notion of $\xi$-monic polynomial,
and we introduce the ring 
$\L_{(\xi)}=
S^{-1}_\xi\L$.
(We take here the occasion to note that for the
first time the localization technique was applied to Novikov rings
and Novikov inequalities in the paper \cite{farber} of M.Farber.
In this paper M.Farber considers the ring $S_\xi^{-1}\L$, where
$\L=\ZZZ[\ZZZ]$, and $\xi$ 
is the inclusion of $\ZZZ$ to $\RRR$.)
Recall next the definition of the Novikov ring $\Lxi$
 (see, e.g. \cite{PaSur}, p. 326).
Denote by
 $\hat {\hat \Lambda}$ 
the abelian group of all the 
linear combinations of the form $\lambda = 
\mathop{\sum}\limits_{g\in \ZZZ^{n+1}} n_g g$  where 
$n_g\in \ZZZ$ and the sum may be infinite. 
Let $\Lxi$ be  the subset of $\hat {\hat \Lambda}$
 consisting of 
$\lambda \in \hat {\hat \Lambda}$ such that for
 every $c\in \RRR$ the set 
$\supp \lambda \cap \xi^{-1}([c,\infty[)$ is finite.
This subset is  called {\it Novikov ring }
(it is not difficult to see that $\Lxi$ has a natural ring structure).

Proceeding to Theorem 3.2, recall that Theorem 1.4 of \cite{PaFree}
asserts that if $\xi$ is injective then $\L_{(\xi)}$ is euclidean.
(The proof is based on a theorem by J.Cl.Sikorav, which asserts that
if $\xi$ is injective then 
$\Lxi$ is euclidean, see \cite{PaFree}, Th. 1.1.)

Therefore we obtain the decomposition
(\ref{Lor}) for 
any monomorphism $\xi$. This implies that (\ref{Lor})
is true for every homomorphism $\eta$
belonging to an open conical set containing $\xi$
(see \cite{PaFree}, the beginning of \S 7).
Since the monomorphisms are dense in $\Hom(\ZZZ^{n+1},\RRR)$,
we obtain the decomposition (\ref{Lor})
for every $\xi$ belonging to
some open and dense conical subset $U$ in 
$\Hom(\ZZZ^{n+1},\RRR)$.

Analyzing further the algebraic structure of the rings 
$\Lxi, \L_{(\xi)}$ (it is done in \cite{PaLor}),
one can prove that $U$ can be chosen in such a way that the complement
$\Hom(\ZZZ^{n+1},\RRR)\sm U$
is a finite union 
$\gR=\cup_i L_i$
of hyperplanes $L_i$. Moreover, each  $L_i$
is of the form
$l_i\otimes\RRR$
where
$l_i\in
\Hom(\ZZZ^{n+1},\ZZZ)$
is an integer hyperplane. 
That proves Theorem 3.2. 
(See \cite{PaLor} for more information 
about the numbers $b_p(C_*,\xi), q_p(C_*,\xi)$.) $\hfill\square$

\pa

The following proposition relates the above numbers and the elements
$a_j^{(p)}$ to the Fitting invariants of 
the boundary operators of
$C_*$. We need some definitions.
Let $A:F_1\to F_2$ be a homomorphism of free finitely generated
$\L$-modules. Let $J_0\subset ... \subset J_r$ be the reduced
Fitting sequence of $A$,  let $\rho_i\in\L$ be the g.c.d. of
the elements of $J_i$, and denote $\rho_i/\rho_{i+1}$ by $\zeta_i(A)$.
Let $\xi:\ZZZ^{n+1}\to\ZZZ$ be a non-zero homomorphism. Denote by
$k(A,\xi)$ the number of those $\rho_i$ which are not $\xi$-monic. Set
$R_j(A,\xi)=\zeta_{k-j}(A,\xi)$ where $k=k(A,\xi)$.
Now let
$C_*=\{0\la ... \la C_{i-1}\lau{\partial_i} C_i\la ... \}$
be a free $\L$-complex.

\prop{
Assume that for
$\xi\in\Hom(\ZZZ^{n+1},\ZZZ)$
and every $p$
 the decomposition
(\ref{Lor})   holds. Then:
1) $b_p(C_*,   \xi)= B_p(C_*)$.\quad
2) $q_p(C_*,   \xi)= k(\partial_{ p+1}   ,\xi)$.\quad
3) For every $p,s$ the elements $a_s^{(p)}$ and
$R_s(\partial_{p+1} ,\xi)$
are equal up to multiplication by a $\xi$-monic element.
4) $Q_p(C_*)$ equals to the number of
not numerically prime
$a_s^{(p)}$.
                    }\label{zeta}

{\it Proof. \quad}  Recall that      the reduced
Fitting sequences
are homotopy invariants of $C_*$. This implies
that 
 $k(\partial_i,\xi)$
and $\k(\partial_i,\xi)$ are  homotopy invariants of $C_*$ for fixed
$\xi$.
 1) is obvious.     Further, let $0\leq p\leq n$ and 
let
$J_0\subset ... \subset J_r$
be the reduced Fitting sequence for $\partial_{p+1}:C_{p+1}\to C_p$.
Then the reduced Fitting sequence $FR(p)$ of
  the localized complex is  a part of
the
sequence
$S_\xi^{-1}J_0\subset ... \subset       S_\xi^{-1}J_r$,
 and the g.c.d. of $S^{-1}_\xi J_i$ is still
$\rho_i$. Therefore the  sequence $FR(p)$
 has $k(\partial_{p+1},\xi)$ terms. Using the principal model
for $C_*$, it is easy to prove that $FR(p)$ equals to the sequence
of principal ideals
$(a_1^{(p)}\cdot ... \cdot a_N^{(p)}),
(a_1^{(p)}\cdot ... \cdot a_{N-1}^{(p)}),
...,
(a_1^{(p)})$ where $N=q_p(C_*,\xi)$. 2), 3)  and 4) follow easily. $\qs$

Let $\xi=l\ov\xi$ where $\ov\xi:\ZZZ^{n+1}\to\ZZZ$
is an epimorphism. Choose an isomorphism
$\Ker\xi\approx\ZZZ^n$, and an element $t\in\ZZZ^{n+1}$ \sut~ $\ov\xi(t)=-1$.
We obtain a decomposition
$\ZZZ^{n+1}=\Ker\xi\oplus\ZZZ$ and
an isomorphism $I(\xi):\L\approx R[t, t^{-1}]$.
Consider the free  $\LL$-complex
${\widehat C}_*(\xi)   =
C_*\tens{\L}\LL$,
where $C_*$ is an $R[t, t^{-1}]$-module via the isomorphism $I(\xi)^{-1}$.
Set
$S(C_*,\xi)=s({\widehat C}_*(\xi))$
( it is easy to check that $S(C_*,\xi)$ depends indeed only 
on $\xi$ and $C_*$).

\theo{
There is
a small        subset $\gM\subset \Hom(\ZZZ^{n+1},\ZZZ)$
\sut~ for every $\xi\notin\gM$ the complex $\wh C_*(\xi)$ is of
asymptotically linear growth and
     $S(C_*,\xi)= B(C_*)+2Q(C_*)$.
         }
{\it Proof.\quad}
The $I(\xi)$-image of  a  $\xi$-monic  polynomial  is  obviously
invertible in $\LL$, therefore
the \ho~
$\L\to R[t, t^{-1}]\to R((t))$ factors
through $\L_{(\xi)}$. Therefore, for every $\xi$ outside
  a small subset $\gN$  the formula (3)
holds with $C_*=\wh C_*(\xi)$.
The complex 
$\wh C_*(\xi)$
is not necessarily of principal type, since the polynomials
$a_s^{(i)}$
in the decomposition (3) 
are not necessarily special. But 
 Proposition \ref{zeta} implies that the elements
$a_s^{(i)}$ in the decomposition (3) can be chosen
between the elements of the finite set
$\{\zeta_j(\partial_{i+1})\}$. Therefore, adding to
$\gN$ some integer hyperplanes if necessary, we can assume that
all $a_s^{(i)}$ are  special. Now  our theorem follows from 2.3.
$\qs$

\section{Proof of the main theorem}

Let $M$ be a closed connected manifold
and $\xi\in H^1(M,\ZZZ)$ be an indivisible cohomology class.
Denote by $\PP_\xi: M(\xi)\to M$ the infinite cyclic covering
\sut~ $\PP_\xi^*(\xi)=0$. Choose a generator 
$t\in\ZZZ\approx\pi_1(M)/\Ker\xi
           $ of the structure
group of $\PP_\xi$
such  that
$\xi(t)=-1$.
Let
$f:M\to S^1$ be a Morse map representing $\xi$, and let
$V$ be its regular
level surface, say, $V=f^{-1}(\l)$. Then $f$ lifts to
a Morse function $F:M(\xi)\to\RRR$ and
     $V$ lifts to $F^{-1}(\l')\subset    M(\xi)$.
    Denote by
$V^-$ the subset $F^{-1}(]-\infty,\l'])$. For $k\geq 1$ denote        by
$W_k$ the cobordism $F^{-1}([\l'-k, \l']),~~
\partial W_k\approx V\sqcup t^kV$,
 and
denote by $\a(k,V)$ its Morse number, that is, the minimal
number of critical points of a Morse function on the cobordism
$W_k$. Note that $\a(k+n, V)\leq \a(k,V)+\a(n,V)$. Therefore
the sequence $\a(k,V)/k$ has a limit as $k\to\infty$.
Denote this limit by $\a(V)$. It is easy to see that
$\a(V)$ depends only on $M$ and $\xi$, so we denote it
by $\a(M,\xi)$.
An elementary construction, using the gluing of the upper part $V$     of
 $\pr W_k$ to the lower part $t^kV$, allows to
obtain the inequality $\MM(M(\xi,k))\leq \a(k,V)+2\MM(V)    $.
In particular, if $\xi$ is represented by a fibration
over $S^1$, the sequence $\MM(M(\xi,k))$ is bounded.

In general, it is
all what we can say about the numbers $\a(k, V)$ and their relation
to the asymptotics of the Morse numbers of cyclic covers.
However, if the fundamental group of $M$ is free abelian
      and $\dim M\geq 6$, one can say
much more.

\prop{
Let $M$ be  a closed connected manifold 
 with a free abelian fundamental group.
Assume that $\dim M\geq 6$.
Let $\xi\in H^1(M),\xi\not=0$. Then the  sequence
$\a(k,V)-      \MM(M(\xi,k))$ is bounded.
         }
\label{varb}

{\it Proof.} Let $\pi_1(M)\approx \ZZZ^{n+1}, n\geq 0   $.
An argument similar to the one used in (\cite{Far}, p.325)
 shows that one can choose
$V$ above such that
 the embedding $V\hookrightarrow M$ induces
an isomorphism $\pi_1(V)\to \Ker\xi$ (such $V$ will be called
{\it  admissible $\xi$-splittings}, see       \cite{PaFree}, p. 371).
In this case all   the embeddings
$V\subset W_k\subset M(\xi)\supset t^kV      $
induce isomorphisms of $\pi_0$ and of $\pi_1$.
Choose an element $T\in\ZZZ^{n+1}$, \sut~ $\xi(T)=-1$.
Let $\overline{M(\xi,k)}=\wi M/T^k$, then
there is  a $\ZZZ^n$-covering 
$\overline{M(\xi,k)}\to M(\xi,k)$.
 Choose a triangulation
of $M$ \sut~ $V$ is a subcomplex of $M$;
then we obtain a $t$-invariant triangulation of
$M(\xi)$ and the corresponding triangulations
of all the covers. There are two exact sequences of
corresponding $\ZZZ[\ZZZ^n]$-complexes:
\begin{align}\label{qqq}
0\ra C_*(\wi V)\ra C_*(\ov{M(\xi,k)})
\ra C_*(\ov{M(\xi,k)},\wi V)\ra 0\\
0\ra C_*(\wi V)\ra C_*(\wi{W_k},\wi{  t^k V})
\ra C_*(\ov{M(\xi,k)},\wi V)\ra 0\label{qqqq}
\end{align}
Proposition  \ref{pryamsum} implies that
there is $C=C(V)$ \sut~ for every $k>0$ we have: \break
$\MM(C_*(\ov{M(\xi,k)}))\geq
\MM(C_*(\wi{W_k},\wi{t^kV})) -C$.
Since
$\MM(C_*(\wi{W_k},\wi{t^kV}))=\a(k, V)$ (see the proof of
 Corollary 6.3 in \cite{Sh1}), our Proposition is proved. $\qs$

{\it Proof of the Main theorem.}
The point 1) follows immediately from \ref{varb}
(with $                           \a(M,\xi)=
\lim_{k\to\infty}\MM(M(\xi,k))/k
$). To prove
2) note that Theorem 3.4
implies that for all $\xi\in H^1(M)$ outside  a small subset
$\gM\subset H^1(M)$ the complex
$(C_*(\wi M))\sphat(\xi)$ is $aslg$, and
$S(C_*(\wi M),\xi)=B(M)+2Q(M)$. Note further that for every admissible
$\xi$-splitting $V$ the complex
$D_*=    C_*(\wi{V^-})\tens{\L} R[[t]]$
 is a basic subcomplex of $C_*(\wi M)\tens{\L} R((t))     $, and that
 $\mu_k(D_*)=      \break
\MM(          C_*(\wi{W_k}, \wi{  t^kV})) =\a(k,V)         $.
Now just apply Proposition 4.1.
$\qs$

\rema{
A similar argument, together with Remark 2.6, shows that for $\xi$
outside  a small subset of $H^1(M)$ the sequence
$(B(M)+2Q(M))k-\MM_\ZZZ(C_*(M(\xi,k)))$ is bounded
                                                   where
$C_*(M(\xi,k))$ is the chain complex of
$M(\xi,k)$, defined over $\ZZZ$ (see 1.2 for the definition of 
$\MM_\ZZZ(\cdot)$).
                  }

\section{Further results and conjectures}

\subsection{Stable Morse numbers}
                    Let $M$ be a closed connected \ma.
Recall  that a {\it stable Morse function} on
 $M$ is
a Morse function $f:M\times \RRR^N\to\RRR$ \sut~
there is a compact $K\subset M\times \RRR^N$,
and a non-degenerate quadratic form $Q$
of index 0
on $\RRR^N$
\sut~ $f(x,y)=Q(y)$ outside $K$.
Let $f:M\times\RRR^N\to\RRR$ be
a stable Morse function. Denote by $\widetilde m_p(f)$ the
number of critical points of $f$
of index $p+N/2$. The
Morse-Pitcher inequalities hold:
$\widetilde m_p(f)\geq b_p(M)     +q_p(M)     +q_{p-1}(M)$.

Denote by $\MM\SS(M)$ the minimal possible number of
critical points of a stable Morse function on $M$; we have
$\MM\SS(M)\leq \MM(M)$.

\theo{
Let $\dim M\geq 6,
\pi_1(M)\approx\ZZZ^{n+1}, n\geq 0
$. There is a subset $\gM\in H^1(M)$
which is  a finite union of integral hyperplanes in $H^1(M)$,
and for every $\xi\notin \gM$ there is
a real number $a$ \sut~ for every $k\in\NNN$ we have
$$
k(B(M)+2Q(M))-a\leq\MM\SS(M(\xi,k))\leq k(B(M)+2Q(M))+a
$$
      }
For the proof just recall that
(by Remark 4.2)
 for every $\xi$ outside  a small subset
of $H^1(M)$
we have
$k(  B(M)+2Q(M))\leq
\MM_\ZZZ(C_*(M(\xi,k))) + C$. $\qs$

We refer to \cite{eligr}
for a systematic exposition of the theory of 
stable Morse functions and 
its applications to Lagrangian intersection theory.

\subsection{ Non generic cohomology classes $\xi\in H^1(M)$ }

Here we construct a manifold $M$ with $\pi_1(M)\approx\ZZZ^2$
and $\dim M\geq 6$, and a class $\xi\in H^1(M)$ \sut~
$\mu(M,\xi)\not= B(M)+2Q(M)$.
Let $N$ be a closed connected \ma~ with $\pi_1(N)\approx\ZZZ,
 \dim N\geq 5$
and  $B(N)\not=0$.
 Set $M=N\times S^1$.
Let $\l:\pi_1(M)\to\ZZZ$,
resp. $\xi:\pi_1(M)\to\ZZZ$, be  epimorphisms
with $\Ker\l=\pi_1(N)$, resp. $\Ker\xi=\pi_1(S^1)$.
Then $\l$ is represented by the fibration $N\times S^1\to S^1$.
Therefore there is an open cone $C\subset H^1(M,\RRR)$
containing $\l$, \sut~ every
integral
non-divisible $\l'\in C$ can be represented
by a fibration, and so $\mu(M,\l')=B(M)+2Q(M)=0$.

Now we shall show that $\mu(M,\xi)\not=0$.
Note that
$M(\xi,k)=N_k\times S^1$ where $N_k$ is the $k$-fold cyclic cover of $N$,
and therefore
$\MM(M(\xi,k))\geq
 \MM_\ZZZ(C_*(N_k\times S^1))$.
We shall obtain a lower estimate for
$\MM_\ZZZ(C_*(N_k\times S^1))$.
Let $\xi_0:\pi_1N\to\ZZZ$ be the restriction
$\xi\mid\pi_1 N$.
Let
$\ov{N}\to N$ be the infinite cyclic covering and
$V\sbs N$ be an admissible $\xi_0$-splitting. Let $W_k$ be the
 corresponding cobordism in $\ov{N}$.
Using exact sequences similar to the exact sequences
(\ref{qqq}, \ref{qqqq}) from  \S 4,
it is easy to prove that
$\MM_\ZZZ(C_*(N_k\times S^1))-
\MM_\ZZZ(C_*(W_k\times S^1,t^k V\times S^1))$
is bounded.
Let $X_k=W_k/t^kV,\quad Y_k=(W_k/t^kV)\times S^1$.
Then
$\MM_\ZZZ(C_*(Y_k))=
\sum_p(b_p(Y_k)+q_p(Y_k)+q_{p-1}(Y_k))$.
Since $H_*(X_k\times S^1)= H_*(X_k)\oplus H_{*-1}(X_k)$
we have
$q_p(X_k\times S^1)\geq q_p(X_k)$, and
$\MM_\ZZZ(C_*(Y_k))\geq \MM_\ZZZ(C_*(X_k))$.
Recall from Proposition 4.1 that the sequence
$\MM(N_k)-\MM_\ZZZ(C_*(X_k))$ is bounded.
Therefore
$\MM_\ZZZ(C_*(N_k\times S^1))
\geq
\MM(N_k)+C\geq
kB(N)+C'$
(where $C$ and $C'$ do not depend on $k$),
and, finally, $\mu(M,\xi)\geq B(N)$.

\subsection{ Non cyclic finite coverings}

The Main Theorem of the present paper allows also
to
deal with some non cyclic finite coverings.

\prop{ Let $M$ be a closed connected manifold, $\dim M\geq 6$,
and $\pi_1(M)\approx \ZZZ^{n+1}, n\geq0$.
Let $M_k\to M$ be the finite covering corresponding to the
subgroup $k\ZZZ^{n+1}$. Then
$$
\lim_{k\to\infty} \frac{\MM(M_k)}{k^{n+1}}= B(M)+2Q(M)
$$
}
{\it Proof. ~~} We shall give only the main idea of the proof.
Let $\xi:\pi_1(M)\to\ZZZ$ be an epimorphism not belonging to
the small set $\gM$ of the Main Theorem.
Then $M_k\to M$ factors through $M(\xi, k)\to M$, and therefore
$\MM(M_k)/k^n\leq k(B(M)+2Q(M))+C$. To obtain the lower estimate,
note that the $\ZZZ^n$-covering
$           \overline{M(\xi,k)}\to M(\xi,k)$ factors through
 $\overline{  M(\xi,k)}\to M_k$ which
 corresponds to the subgroup $G_k=k\ZZZ^n\subset \ZZZ^n$.
Therefore
$\MM(M_k)\geq \MM_{\ZZZ[G_k]}(C_*(\overline{M(\xi,k)}))$.
To obtain the lower estimate for
$\MM_{\ZZZ[G_k]}(C_*(\overline{M(\xi,k)}))$,  use (7) and (8)
to reduce the question to finding the corresponding lower
estimate
 for
$\MM_{\ZZZ[G_k]} (C_*(\wi W_k, \wi{ t^k    V}))$.
Then
proceed similarly to the proof of the Main Theorem. $\qs$

It seems that a similar result
 must hold for more general
systems of non-cyclic
finite coverings.
To discuss a more general setting we need some definitions.
Let $G$ be a  group. A sequence
of subgroups
$G=G_0\supset G_1\supset ...$
will be called  a {\it tower} if
for every  $i$ the index of $G_i$ in $G$ is finite.
 It will be called a
{\it nested tower} if, moreover,  $\cap_nG_n=\{0\}$.

If $M$ is a
closed connected
\ma~ with $\pi_1(M)=G$ and
 $\gg=     \{G_n\}$ is a  tower of
subgroups of
 $\pi_1(M)$, then
there is
the corresponding
 tower of finite coverings
$M = M_0 \la ...\la M_n ...  $~~ of $M$.
The sequence $\MM(M_k)/\mid G/G_k \mid$ is decreasing, \th~ it has a limit
which will be denoted by $\mu(\gg)$.
Recall a theorem of W.L\"uck \cite{Luck}, saying
that if $\gg$ is nested,
then
 the limit of the sequence
$b_p(M_k)/\mid G/G_k\mid    $ exists
and is equal to $b_p^{(2)}(M)$.
\vskip0.1in
{\bf Problem.}\hspace{0,2cm} Is it true in general (at least 
for $\dim M\geq 6$)
that $\mu(\gg)$ does not depend on the choice of
the nested tower $\gg$?
\vskip0.1in
We believe that $\mu(\gg)$ does not depend on $\gg$ for the case of free
 abelian fundamental group. Here is a result in this direction.
Let $G=\ZZZ^{n+1}$. For a tower
$\gg=\{G=G_0\supset G_1\supset ...\}$
we denote
$\max_i m(G/G_i)$ by $r(\gg)$ (here $m(H)$ stands for the minimal
number of generators of $H$).  Denote by $\gg^{[i]}$ the tower
$G_i\supset G_{i+1}\supset ...$
The sequence $r(\gg^{[i]})$ is decreasing, denote its limit
by $\rho(\gg)$; then $\rho(\gg)\leq \rk G$.

\prop{
Let $M$ be a closed connected manifold, $\dim M\geq 6$,
$\pi_1(M)\approx \ZZZ^{n+1}, n\geq 0$.
Let $\gg=\{G_k\}$ be a tower, and let $M_k\to M$ be the finite covering
corresponding to $G_k\subset \pi_1(M)$.
Assume that                   $\rho(\gg)=n+1$.
Then
\begin{equation}\label{predel}
\lim_{k\to\infty} \frac{\MM(M_k)}{\mid G/G_k \mid}= B(M)+2Q(M)
\end{equation}
              }
{\it Proof. } It is not difficult to show that if $\rho(\gg)=n+1$,
 then for each $k$
there are $a_k, b_k\in\NNN$ \sut
 $a_kG\supset G_k\supset b_kG$
with $a_k, b_k\to\infty$ as $k\to \infty$.
 Now our Proposition follows from Proposition 5.2.
$\qs$
\vskip0.1in
{\bf Conjecture.~~}
Equality
(\ref{predel}) is true for every
nested tower in $\ZZZ^{n+1}$.

\end{document}